\documentclass[11pt,reqno]{amsart}

\usepackage[T1]{fontenc}
\usepackage{lmodern,microtype}
\usepackage[margin=1.03in]{geometry}
\usepackage{amsmath,amssymb,amsthm}
\usepackage{xcolor}
\usepackage[colorlinks=true,linkcolor=blue!50!black,
  citecolor=blue!50!black,urlcolor=blue!50!black]{hyperref}

\newtheorem{theorem}{Theorem}[section]
\newtheorem{lemma}[theorem]{Lemma}
\newtheorem{proposition}[theorem]{Proposition}
\newtheorem{corollary}[theorem]{Corollary}
\DeclareMathOperator{\Ann}{Ann}
\DeclareMathOperator{\End}{End}
\DeclareMathOperator{\Frac}{Frac}
\DeclareMathOperator{\Hom}{Hom}
\newcommand{\F}{\mathbf F}
\newcommand{\m}{\mathfrak m}
\newcommand{\n}{\mathfrak n}

\title[Divisorial domains and normalization]{A divisorial domain with non-Pr\"ufer normalization}

\author{Viet-Hoang Tran}
\address{Department of Mathematics, National University of Singapore, Singapore 119076}
\email{hoang.tranviet@u.nus.edu}
\urladdr{https://vh-tran.github.io/}
\date{}
\hypersetup{pdftitle={A divisorial domain with non-Prüfer normalization},
  pdfauthor={Viet-Hoang Tran}}

\begin{document}
\raggedbottom
\begin{abstract}
We answer a 1968 question of Heinzer by constructing a two-dimensional local
divisorial domain whose normalization is not Pr\"ufer.
\end{abstract}
\maketitle

\section{Heinzer's question and its status in the literature}

Let \(D\) be a domain with fraction field \(K\). For \(D\)-submodules \(U,V\)
of \(K\), write \((U:V)=\{z\in K\mid zV\subseteq U\}\). The divisorial closure
of a nonzero fractional ideal \(I\) is \(I_v=(D:(D:I))\). The domain \(D\)
is \emph{divisorial} if \(I=I_v\) for every such \(I\). Its normalization
\(\overline D\) is its integral closure in \(K\). In 1968 Heinzer \cite{Heinzer}
asked whether \(\overline D\) must be a Pr\"ufer domain whenever \(D\)
is divisorial. A Pr\"ufer domain is one in which each nonzero finitely generated
ideal is invertible.

Two classical cases give an affirmative answer. An integrally closed
divisorial domain is Pr\"ufer \cite[Theorem~5.1]{Heinzer}. A Noetherian
divisorial domain that is not a field has dimension one
\cite[Corollary~4.3]{Heinzer}, so its normalization is Dedekind by the
Krull--Akizuki theorem. Gabelli, Houston and Picozza
\cite[Questions~3.7 and~3.8]{GHP} restated Heinzer's question and proved it
equivalent to asking whether every \(w\)-divisorial domain is a UM\(t\)-domain.
A domain is \(w\)-divisorial if each of its \(w\)-ideals is divisorial.
In particular, every divisorial domain is \(w\)-divisorial.

In 2026 Tran, Vo and Nguyen \cite{TVN} answered a question of Vasconcelos
by constructing a one-dimensional coherent local domain with non-Pr\"ufer
normalization. Kim
\cite[Theorem~A]{Kim} proved that the same domain is Gorenstein Pr\"ufer.
In particular, each of its nonzero finitely generated ideals is divisorial.
Heinzer's question requires the identity \(I=I_v\) also for infinitely
generated ideals.

\begin{theorem}\label{thm:main}
For each prime \(p\), put \(k=\F_p(t_1,t_2,\ldots)\), where the \(t_i\)
are algebraically independent. There exists a two-dimensional local divisorial
domain \(R\) of characteristic \(p\) with normalization \(S=k[[x,y]]\).
In particular, \(\overline R\) is not Pr\"ufer, and \(R\) fails the
UM\(t\) condition.
\end{theorem}

The construction uses a derivation to obtain a subring of a prescribed
normal domain, an approach also used by Olberding \cite{Olberding}.
Here the derivation takes values in the injective hull of the residue field
of \(S\). This choice controls annihilators in every quotient \(R/cR\)
by a nonzero nonunit.

\section{A derivation with prescribed values}

Write \(\n=(x,y)S\) and let \(E=E_S(k)\) be the injective hull of \(k=S/\n\).
For an ideal \(H\) and a module \(M\) over the same ring, put
\(M[H]=\{u\in M\mid Hu=0\}\), and write \(M[c]=M[(c)]\).

\begin{lemma}\label{lem:derivation}
There is an \(\F_p\)-derivation \(\delta:S\to E\) which is surjective as a map
of additive groups.
\end{lemma}

\begin{proof}
The equality \(E=\bigcup_{j\ge1}E[\n^j]\) holds over the Noetherian local ring
\(S\) \cite[Lemma~47.7.3]{Stacks}. Each
\(E[\n^j]=\Hom_S(S/\n^j,E)\) has finite length. Indeed \(\Hom_S(-,E)\) is
exact and sends the simple module \(k\) to \(k\), the socle of \(E\).
As \(k\) is countable, so is \(E\). Enumerate its elements as \(e_1,e_2,\ldots\).

In the ordinary module of K\"ahler differentials \(\Omega_{S/\F_p}\), the
\(dt_i\) are \(S\)-linearly independent. For each \(j\), the derivation
\(\partial/\partial t_j\) of \(k\) extends coefficientwise to \(S\):
the product rule holds since each coefficient of a product is a finite sum.
The corresponding \(S\)-linear functional on \(\Omega_{S/\F_p}\) takes \(dt_i\)
to the Kronecker delta \(\delta_{ij}\). On the free submodule
\(\bigoplus_i S\,dt_i\), send \(dt_i\) to \(e_i\). Injectivity of \(E\)
extends this map to all of \(\Omega_{S/\F_p}\). Its composite with the
universal derivation has \(\delta(t_i)=e_i\), so its image contains every
element of \(E\).
\end{proof}

Fix such a \(\delta\) and put \(R=\ker\delta\). This is a local domain with
maximal ideal \(\m=R\cap\n\), since an element \(r\in R\) invertible in \(S\)
satisfies \(\delta(r^{-1})=-r^{-2}\delta(r)=0\).
In characteristic \(p\) we have \(S^p\subseteq R\), so \(S\) is integral
over \(R\). For each \(s\in S\), a power of \(\n\) kills \(\delta(s)\).
Choose \(N\) with \(x^{pN}\delta(s)=0\). Then \(x^{pN}\) and \(x^{pN}s\)
both belong to \(R\), so \(\Frac R=\Frac S\). Since \(S\) is integrally
closed, \(\overline R=S\). Integrality gives \(\dim R=\dim S=2\).
In a local Pr\"ufer domain every nonzero finitely generated ideal is
principal, whereas \(\dim_k(\n/\n^2)=2\). Thus \(S\) is not Pr\"ufer.

\section{Divisoriality of all ideals}

\begin{lemma}\label{lem:ann}
Let \(J\) be a square-zero ideal of a commutative ring \(B\).
Suppose \(A=B/J\) is complete Noetherian local with residue field \(\ell\)
and \(J\simeq E_A(\ell)\) as an \(A\)-module. Then \(B\simeq E_B(\ell)\)
as \(B\)-modules. Every ideal \(L\) of \(B\) satisfies
\(L=\Ann_B(\Ann_B L)\).
\end{lemma}

\begin{proof}
Since \(J\) is nilpotent, \(B\) is local with residue field \(\ell\).
Set \(W=E_B(\ell)\). Its submodule \(W[J]=\Hom_B(A,W)\) is injective over
\(A\) by adjunction and contains \(\ell\) essentially. Hence it is
\(E_A(\ell)\). Fix an isomorphism \(\theta:J\to W[J]\) and extend it by
injectivity to a \(B\)-linear map \(f:B\to W\).

Matlis duality identifies scalar multiplication \(A\to\End_A(J)\) with
the completion map \cite[Lemma~47.7.5]{Stacks}, hence with an isomorphism.
In particular \(J\) is faithful over \(A\). If \(f(b)=0\), then
\(\theta(bj)=jf(b)=0\) for every \(j\in J\). Faithfulness gives \(b\in J\),
where \(f\) is injective. To see surjectivity, take \(w\in W\).
The map \(j\mapsto jw\) has image in \(W[J]\) because \(J^2=0\).
Composing with \(\theta^{-1}\) gives an endomorphism of \(J\), hence
multiplication by some \(\bar b\in A\). Lift it to \(b\in B\). Then
\(jw=\theta(bj)=jf(b)\) for every \(j\), so \(w-f(b)\in W[J]=f(J)\).

Finally, if \(b\notin L\), the nonzero cyclic submodule generated by \(b+L\)
in \(B/L\) has a quotient \(\ell\) in which its generator survives. Embed
\(\ell\) in \(W\) and extend the resulting homomorphism to \(B/L\).
Under \(W\simeq B\), the image of \(1+L\) annihilates \(L\) and has
nonzero product with \(b\). This proves the double-annihilator identity.
\end{proof}

\begin{proposition}\label{prop:quotient}
For \(0\ne c\in\m\), the map \(B=R/cR\to A=S/cS\) is onto. Its kernel
\(J\) is square-zero and is isomorphic to \(E_A(k)\) as an \(A\)-module.
\end{proposition}

\begin{proof}
Multiplication by \(c\) on \(E\) is onto. For any \(e\in E\), the
\(S\)-linear map \(cS\to E\) sending \(c\) to \(e\) extends to \(S\),
as \(c\) is regular and \(E\) is injective. Lemma~\ref{lem:derivation}
therefore gives, for any \(s\in S\), some \(t\in S\) with
\(c\delta(t)=\delta(s)\). Then \(s-ct\in R\), proving surjectivity.

The kernel is \(J=(cS\cap R)/cR\). If \(cs,ct\in R\), then
\(c\delta(s)=c\delta(t)=0\), whence \(\delta(cst)=0\) and
\((cs)(ct)\in cR\). Thus \(J^2=0\). The assignment
\[
     J\longrightarrow E[c],\qquad cs+cR\longmapsto\delta(s)
\]
is well-defined and injective because \(\ker\delta=R\) and \(c\) is regular.
It is onto because \(\delta(S)=E\), and \(cs\in R\) whenever
\(c\delta(s)=0\). It is \(A\)-linear since \(R\to A\) is onto and
\(\delta(rs)=r\delta(s)\) for \(r\in R\).
Finally \(E[c]=\Hom_S(A,E)\) is injective over \(A\) and contains \(k\)
essentially, so it is \(E_A(k)\).
\end{proof}

Let \(I\) be any nonzero proper ideal of \(R\), and choose \(0\ne c\in I\).
The ring \(S/cS\) is complete Noetherian local, so
Lemma~\ref{lem:ann} applies to \(B=R/cR\) by
Proposition~\ref{prop:quotient}. For \(cR\subseteq H\subseteq R\) the
fractional colon \((cR:H)\) lies in \(R\), because \(c\in H\), and
contains \(cR\). Consequently
\[
    \Ann_B(H/cR)=(cR:H)/cR.
\]
Applying the double-annihilator identity to \(I/cR\) yields
\[
    I=(cR:(cR:I))=(cR:c(R:I))=(R:(R:I))=I_v .
\]
The identity is immediate for \(I=R\). Clearing a denominator extends it
to every nonzero fractional ideal because \((aI)_v=aI_v\) for
\(a\in\Frac(R)^\times\). This proves the divisoriality assertion of
Theorem~\ref{thm:main}.

\section{An upper to zero}

A nonzero prime of \(D[Z]\) contracting to zero in \(D\) is called an
\emph{upper to zero}. A domain \(D\) is a UM\(t\)-domain if each upper to
zero is a maximal \(t\)-ideal of \(D[Z]\). Here \(H_t\) is the union of
\(F_v\) over the nonzero finitely generated subideals \(F\subseteq H\).
In particular every nonzero divisorial ideal is a \(t\)-ideal, since
\(H\subseteq H_t\subseteq H_v\).

\begin{corollary}\label{cor:upper}
In \(R[Z]\), the kernel \(Q\) of evaluation at \(y/x\in\Frac R\) is an
upper to zero strictly contained in the proper \(t\)-ideal \(\m[Z]\).
\end{corollary}

\begin{proof}
Evaluation is injective on \(R\) and its target is a field. Its kernel is
nonzero since \(x^pZ^p-y^p\in Q\). If \(\sum_{i=0}^{d}r_iZ^i\in Q\), then
\(\sum_{i=0}^{d}r_i x^{d-i}y^i=0\) in \(S\). The monomials
\(x^{d-i}y^i\) give a \(k\)-basis of \(\n^d/\n^{d+1}\), so the residue of
every \(r_i\) in \(k\) is zero. Thus \(Q\subseteq\m[Z]\), and \(x^p\)
witnesses strict inclusion.

Since \(\m=\m_v\), it is an intersection of principal fractional ideals
of \(R\). Extending each to \(R[Z]\) preserves this intersection, as
membership is determined coefficientwise. Hence \(\m[Z]\) is divisorial
in \(R[Z]\), and therefore a proper \(t\)-ideal. This proves the remaining
assertion of Theorem~\ref{thm:main}.
\end{proof}

\section*{Acknowledgment}
The counterexample and its proof were generated through discussions with
ChatGPT~6 Astra, which also assisted in drafting the manuscript.

\bibliographystyle{amsplain}
\begingroup
\raggedright
\bibliography{references}
\endgroup
\end{document}